\documentclass[12pt,reqno]{amsart}
\usepackage{amsmath}
\usepackage{amssymb}
\usepackage[left=1in,top=1in,right=1in,bottom=1in]{geometry}
\usepackage[usenames]{color}
\usepackage{enumerate}
\usepackage[normalem]{ulem}
\usepackage{soul}

\begin{document}
\newtheorem{theorem}{Theorem}[section]
\newtheorem{lemma}[theorem]{Lemma}
\newtheorem{definition}[theorem]{Definition}
\newtheorem{conjecture}[theorem]{Conjecture}
\newtheorem{proposition}[theorem]{Proposition}
\newtheorem{claim}[theorem]{Claim}
\newtheorem{algorithm}[theorem]{Algorithm}
\newtheorem{corollary}[theorem]{Corollary}
\newtheorem{observation}[theorem]{Observation}
\newtheorem{problem}[theorem]{Open Problem}
\newtheorem{question}[theorem]{Question}
\newcommand{\noin}{\noindent}
\newcommand{\ind}{\indent}
\newcommand{\om}{\omega}
\newcommand{\pp}{\mathcal P}
\newcommand{\AC}{\mathcal A \mathcal C}
\newcommand{\bAC}{\overline{\AC}}
\newcommand{\ppp}{\mathfrak P}
\newcommand{\N}{{\mathbb N}}
\newcommand{\Z}{{\mathbb Z}}
\newcommand{\LL}{\mathbb{L}}
\newcommand{\R}{{\mathbb R}}
\newcommand{\E}{\mathbb E}
\newcommand{\Prob}{\mathbb{P}}
\newcommand{\eps}{\varepsilon}
\newcommand{\ram}[1]{\hat{R}({#1})}
\newcommand{\G}{{\mathcal{G}}}
\newcommand{\bin}{{\mathrm{Bin}}}

\newcommand{\Ss}{{\mathcal S}}
\newcommand{\Nn}{{\mathcal N}}

\newcommand{\ceil}[1]{\left \lceil #1 \right \rceil}
\newcommand{\floor}[1]{\left \lfloor #1 \right \rfloor}
\newcommand{\size}[1]{\left \vert #1 \right \vert}
\newcommand{\dist}{\mathrm{dist}}

\title{Note on the multicolour size-Ramsey number for paths}

\author{Andrzej Dudek}
\address{Department of Mathematics, Western Michigan University, Kalamazoo, MI, USA}
\thanks{The first author was supported in part by Simons Foundation Grant \#522400.}
\email{\tt andrzej.dudek@wmich.edu}

\author{Pawe\l{} Pra\l{}at}
\address{Department of Mathematics, Ryerson University, Toronto, ON, Canada}
\thanks{The second author was supported in part by NSERC}
\email{\tt pralat@ryerson.ca}

\begin{abstract}
The size-Ramsey number $\ram{F,r}$ of a graph $F$ is the smallest integer $m$ such that there exists a graph $G$ on $m$ edges with the property that any colouring of the edges of $G$ with $r$ colours yields a monochromatic copy of $F$. In this short note, we give an alternative proof of the recent result of Krivelevich that $\ram{P_n,r} = O((\log r)r^2 n)$. This upper bound is nearly optimal, since it is also known that $\ram{P_n,r} = \Omega(r^2 n)$.
\end{abstract}

\maketitle

\section{Introduction}
Following standard notation, we write $G\to (F)_r$ if any $r$-edge colouring of $G$ (that is, any colouring of the edges of $G$ with $r$ colours) yields a monochromatic copy of $F$. We define the \emph{size-Ramsey number} of $F$ as $\ram{F,r} = \min \{ |E(G)| : G \to (F)_r \}$; that is, $\ram{F,r}$ is the smallest integer $m$ such that there exists a graph $G$ on $m$ edges such that $G \to (F)_r$. For two colours (that is, for $r=2$) the size-Ramsey number was first studied by Erd\H os, Faudree, Rousseau and Schelp~\cite{EFRS1978}.

In this note, we are concerned with the size-Ramsey number of the path $P_n$ on $n$ vertices. It is obvious that $\ram{P_n,2} = \Omega(n)$ and it is easy to see that $\ram{P_n,2} = O(n^2)$; for example, $K_{2n}\to (P_n)_2$. 
The exact behaviour of $\ram{P_n,2}$ was not known for a long time. In fact, Erd\H{o}s~\cite{E81} offered \$100 for a proof or disproof that $\ram{P_n,2} / n \to \infty$ and $\ram{P_n,2} / n^2 \to 0$. This problem was solved by Beck~\cite{B83} in 1983 who, quite surprisingly, showed that $\ram{P_n,2} < 900 n$. (Each time we refer to inequality such as this one, we mean that the inequality holds for sufficiently large $n$.) A variant of his proof, provided by Bollob\'{a}s~\cite{B01}, gives $\ram{P_n,2} < 720 n$. Recently, the authors of this paper~\cite{DP15} used a different and more elementary argument that shows that $\ram{P_n,2} < 137 n$. The argument was subsequently tuned by Letzter~\cite{L15} who showed that $\ram{P_n,2} < 91 n$, and then further refined by the authors of this paper~\cite{DP2017} who showed that $\ram{P_n,2} \le 74n$.
On the other hand, the first nontrivial lower bound was provided by Beck~\cite{B90} and his result was subsequently improved by Bollob\'as~\cite{B86} who showed that $\ram{P_n,2} \ge (1+\sqrt{2}) n - O(1)$. The strongest lower bound, $\ram{P_n,2} \ge 5n/2 - O(1)$, was proved in~\cite{DP2017}.

Let us now move to the multicolour version of this graph parameter. It was proved in~\cite{DP2017} that $\frac{(r+3)r}{4}n-O(r^2) \le \ram{P_n,r} \le 33 r 4^r n$. It follows that $\ram{P_n,r}$ is linear for any fixed value of $r$ but the two bounds are quite apart from each other in terms of their dependence on $r$. Subsequently, Krivelevich~\cite{K2016} showed that in fact the dependence on $r$ is (nearly) quadratic; that is, $\ram{P_n,r} = r^{2+o_r(1)}n$. Here is the precise statement of his result: 
\begin{theorem}[\cite{K2016}]
For any $C>5$, $r\ge 2$, and all sufficiently large $n$ we have
\[
\ram{P_n,r} < 400^5 C r^{2+\frac{1}{C-4}}n.
\]
\end{theorem}
\noindent
It is straightforward to see that $C=C(r)$ that minimizes the upper bound in this theorem is of order $\log r$. As a result we get that $\ram{P_n,r} = O((\log r)r^2 n)$. In this note, we give an alternative proof of this fact.
\begin{theorem}\label{thm:main}
For any integer $r\ge 2$ and all sufficiently large $n$ we have
\[
\ram{P_n,r} < 600 (\log r) r^2 n.
\]
\end{theorem}
\noindent
It will follow from the proof that the constant 600 is not optimal. Since we believe that the factor $\log r$ is not necessary, we do not attempt to optimize it.

\section{Proof}\label{sec:proof}

Before we move to the proof of Theorem~\ref{thm:main}, we need one, straightforward, auxiliary result.
\begin{proposition}\label{prop:G}
For any integer $r\ge 2$ there exists an integer $N=N(r)$ such that the following holds. For any integer $n \ge N$, there exists a graph $G=(V,E)$ such that 
\begin{itemize}
\item [(i)] $|V| = 7rn$, 
\item [(ii)] $500 (\log r) r^2 n < |E| < 600 (\log r) r^2 n$, and 
\item [(iii)] for every two disjoint sets $S,T \subseteq V$, $|S|=|T|=n$, the number of edges induced by $S\cup T$ with at least one endpoint in $S$ is at most $70(\log r)n$.
\end{itemize}
\end{proposition}

\begin{proof}
The proof is an easy application of random graphs.
Recall that the \emph{binomial random graph} $\G(n,p)$ is a distribution over the class of graphs with vertex set $[n]$ in which every pair $\{i,j\} \in \binom{[n]}{2}$ appears independently as an edge in $G$ with probability~$p$, which may (and usually does) tend to zero as $n$ tends to infinity. Furthermore,  
we say that events $A_n$ in a probability space hold \emph{asymptotically almost surely} (or \emph{a.a.s.}), if the probability that $A_n$ holds tends to $1$ as $n$ goes to infinity. 

Fix any integer $r\ge 2$. It suffices to show that the random graph $G\in \G(7rn, p)$ with $p = 22(\log r)/n$ a.a.s.\ satisfies properties (ii) and (iii). (Property (i) trivially holds.)  Indeed, if this is the case, then there exists an integer $N=N(r)$ such that the desired properties hold with probability at least $1/2$ for $G\in \G(7rn, p)$ for all $n \ge N$. This implies that for each $n \ge N$, there exists at least one graph with these properties. 

\medskip \emph{Property (iii)}: Fix any two disjoint subsets $S, T \subseteq V$, both of cardinality $n$. Let $X_{S,T}$ be the random variable counting the number of edges induced by $S\cup T$ with at least one endpoint in $S$. Clearly, $X_{S,T}$ has the binomial distribution $\bin\big{(} |S|\cdot |T| + \binom{|S|}{2}, p \big{)}$ with $\E(X_{S,T}) = (3/2+o(1)) n^2p = (33+o(1)) (\log r)n$. It follows from Chernoff's bound (see, for example, Corollary 21.7 in~\cite{FK16}) that 
\[
\Pr(X_{S,T} \ge 70(\log r)n) \le \Pr(X_{S,T} \ge 2 \E(X_{S,T})) \le \exp(-\E(X_{S,T})/3) \le \exp(-10.9(\log r)n).
\]
Thus, the probability that there exist $S$ and $T$ such that $X_{S,T} \ge 70(\log r)n$ is, by the union bound, at most
\begin{eqnarray*}
\binom{7rn}{n}^2 \exp(-10.9(\log r) n) &\le& (7er)^{2n} \exp(-10.9(\log r) n) \\
&\le& \exp \left( n \Big( 2\log(7e r) - 10.9 \log r \Big) \right) = o(1),
\end{eqnarray*}
since $(7e r)^2 < r^{10.9}$ for any $r\ge 2$. Property (iii) holds a.a.s.

\medskip \emph{Property (ii)}: This property is straightforward to prove. Note that $|E|$ is distributed as $\bin\big{(} \binom{7rn}{2}, p \big{)}$ with $\E(|E|) = (539+o(1)) (\log r) r^2 n$. It follows immediately from Chernoff's bound that property (ii) holds a.a.s. The proof of the proposition is finished. 
\end{proof}

\begin{proof}[Proof of Theorem~\ref{thm:main}]
The proof is based on the depth first search algorithm (DFS), applied several times, and it is a variant of the previous approach taken in~\cite{DP2017} where it was proved that $ \ram{P_n,r} \le 33 r 4^r n$. Using the DFS algorithms in a Ramsey-type problem was first successfully applied by Ben-Eliezer, Krivelevich and Sudakov~\cite{BKS2012}.

Fix $r\ge 2$ and suppose that $n$ is sufficiently large so that Proposition~\ref{prop:G} can be applied. Let $G=(V,E)$ be a graph satisfying properties (i)--(iii) from Proposition~\ref{prop:G}. We will show that $G\to (P_n)_r$ that implies the desired upper bound as $|E| < 600 (\log r) r^2 n$ by property (ii). Consider any $r$-colouring of the edges of $G$. By averaging argument, there is a colour (say blue) such that the number of blue edges is at least $|E(G)|/r$. For a contradiction, suppose that there is no monochromatic copy of $P_n$; in particular, there is no blue copy of~$P_n$. 

From now on, we restrict ourselves to the graph $G_1 = (V_1=V, E_1 \subseteq E)$, the subgraph of $G$ induced by blue edges. We perform the following algorithm on $G_1$ to construct a path $P$. Let $v_1$ be an arbitrary vertex of $G_1$, let $P=(v_1)$, $U = V \setminus \{v_1\}$, and $W = \emptyset$. If there exists an edge from $v_1$ to some vertex in $U$ (say from $v_1$ to $v_2$), we extend the path as $P=(v_1,v_2)$ and remove $v_2$ from $U$. We continue extending the path $P$ this way for as long as possible. Since there is no $P_n$ in the blue graph, we must reach a point of the process in which $P$ cannot be extended, that is, there is a path from $v_1$ to $v_k$ ($k < n$) and there is no edge from $v_k$ to $U$ (including the case when $U$ is empty). This time, $v_k$ is moved to $W$ and we try to continue extending the path from $v_{k-1}$, reaching another critical point in which another vertex will be moved to $W$, etc. If $P$ is reduced to a single vertex $v_1$ and no edge to $U$ is found, we move $v_1$ to $W$ and simply re-start the process from another vertex from $U$, again arbitrarily chosen.

Observe that during this algorithm there is never an edge between $U$ and $W$. Moreover, in each step of the process, the size of $U$ decreases by 1 or the size of $W$ increases by 1. The algorithm ends when $U$ becomes empty and all vertices from $P$ are moved to $W$. However, we will finish it prematurely, distinguishing $7r$ phases; phase $i$ starts with graph $G_i=(V_i,E_i)$ and ends when for the first time $|W| = n$. Before we move to the next phase, we set $S_i = W$, $T_i = V(P)$, and $F_i$ to be all edges incident to $W$. Then, we set $V_{i+1} = V_i \setminus W$ and $G_{i+1} = G_i [V_{i+1}]$, the graph induced by $V_{i+1}$ (in other words, $G_{i+1}$ is formed from $G_i$ by removing vertices from $W$ together with $F_i$, all edges incident to them). Phase $i$ ends now and we move to phase $i+1$ where we run the algorithm on $G_{i+1}$. 

There are a few important observations. Note that, by property (i), $|V|=|V_1|=7rn$ so the last phase, phase $7r$, finishes with $U=\emptyset$ and $T_{7r}=\emptyset$. As a result, family $(F_i : 1 \le i \le 7r)$ is a partition of $E_1$. By construction, $|S_i| = n$ for all $i$ and, since there is no path on $n$ vertices in $G_1$ (and so also in any $G_i$), $|T_i| < n$ for all $i$.  Hence, $|F_i| < 70(\log r)n$ by property (iii). Putting these things together and using property (ii) in the very last inequality, we get the desired contradiction:
\begin{eqnarray*}
|E|/r \le  |E_1| = |F_1| + |F_2| + \dots + |F_{7r}| \le 7r \cdot 70(\log r)n < 500 (\log r) r n < |E|/r.
\end{eqnarray*}
The proof is finished.
\end{proof} 

\bibliographystyle{amsplain}

\providecommand{\bysame}{\leavevmode\hbox to3em{\hrulefill}\thinspace}
\providecommand{\MR}{\relax\ifhmode\unskip\space\fi MR }
% \MRhref is called by the amsart/book/proc definition of \MR.
\providecommand{\MRhref}[2]{%
  \href{http://www.ams.org/mathscinet-getitem?mr=#1}{#2}
}
\providecommand{\href}[2]{#2}

\end{document}